\input{epsf}

\newcommand{\et}{\hspace{-0.08in}{\bf .}\hspace{0.1in}}
\newcommand{\BOX}{\hbox {$\sqcap$ \kern -1em $\sqcup$}}

\newcounter{letter}     
\newenvironment{alphalist}{\begin{list}
{{\normalshape(\alph{letter})}}{\usecounter{letter}}}{\end{list}}               

\newcommand {\T}{\cal T}
\newcommand {\C}{\cal C}
\newcommand {\D}{\cal D}

\newcommand {\tensor}{\otimes}

\newcommand {\hcomp}{\circ}
\newcommand {\btensor}{\bigotimes}
\newcommand {\tto}{\Rightarrow}
\newcommand {\maps}{\colon}
\newcommand {\from}{\colon}
\newcommand {\R}{{\bf R}}

\newtheorem{thm}{Theorem}    

\newtheorem{defn}[thm]{Definition}

        \newcommand{\be}{\begin{equation}}
        \newcommand{\ee}{\end{equation}}
        \newcommand{\ba}{\begin{eqnarray}}
        \newcommand{\ea}{\end{eqnarray}}
        \newcommand{\ban}{\begin{eqnarray*}}
        \newcommand{\ean}{\end{eqnarray*}}
        \newcommand{\barr}{\begin{array}}
        \newcommand{\earr}{\end{array}}

\documentstyle[12pt]{article}
\begin{document}

      \begin{center}
      {\bf 2-Tangles \\}
      \vspace{0.5cm}
      {\em John C.\ Baez and Laurel Langford\\}
      \vspace{0.3cm}
      {\small Department of Mathematics,  University of California\\ 
      Riverside, California 92521 \\
      USA\\ }
      \vspace{0.3cm}
      {\small email: baez@math.ucr.edu, langford@math.ucr.edu\\}
      \vspace{0.3cm}
      {\small March 17, 1997 \\ }
      \end{center}

\begin{abstract}
Just as links may be algebraically described as certain morphisms in
the category of tangles, compact surfaces smoothly embedded in $\R^4$
may be described as certain 2-morphisms in the 2-category of
`2-tangles in 4 dimensions'.  In this announcement we give a purely
algebraic characterization of the 2-category of unframed unoriented
2-tangles in 4 dimensions as the `free semistrict braided monoidal
2-category with duals on one unframed self-dual object'.  A
forthcoming paper will contain a proof of this result using the movie
moves of Carter, Rieger and Saito.  We comment on how one might use
this result to construct invariants of 2-tangles.
\end{abstract}

\section{Introduction}

Recent work on `quantum invariants' of knots, links, tangles, and
3-manifolds depends crucially on a purely algebraic characterization
of tangles in 3-dimensional space.  It follows from work of Freyd
and Yetter, Turaev, and Shum \cite{FY,Shum,T,Y} that isotopy classes of 
framed oriented tangles in 3 dimensions are the morphisms of a certain 
very special category: the `free braided monoidal category with duals on one
object'.  It follows that we can easily obtain functors from this
category to other braided monoidal categories with duals, such as the
category of representations of a quantum group.  Any such functor
gives an invariant of tangles, and therefore of knots and links.

The `tangle hypothesis' \cite{BD} suggests a vast generalization of
this result, applicable to $n$-manifolds smoothly embedded in
$(n+k)$-dimensional space.  It says that framed oriented $n$-tangles
in $n+k$ dimensions are the $n$-morphisms of the `free $k$-tuply
monoidal weak $n$-category with duals on one object'.  The hope is
that a precise formulation and proof of this hypothesis will open the
door to constructing quantum invariants of $n$-dimensional submanifolds
of $\R^{n+k}$.  

Unfortunately the tangle hypothesis involves concepts from topology
and $n$-category theory that have so far only been worked out in
certain low-dimensional cases.  As a kind of warmup, we wish to prove
a version of this hypothesis in the case $n = k = 2$.  So far we have
only completed work on the unframed, unoriented case, which allows us
to take maximal advantage of the recent work of Carter, Rieger and
Saito \cite{CRS}.  Since the theory of $k$-tuply monoidal weak
$n$-categories is not yet well developed for $n = k = 2$, we use the
better-understood `semistrict' ones as a kind of stopgap.  These are
also known as `semistrict braided monoidal 2-categories'.  The result
announced here is thus that the 2-category of unframed unoriented
2-tangles in 4 dimensions is the `free semistrict braided monoidal
2-category with duals on one unframed self-dual object'.  The proof
appears in Langford's dissertation \cite{L}, and will be published as
part of the `Higher-Dimensional Algebra' series \cite{BL}.

To appreciate this result, one needs some feeling for the topology and
algebra involved: that is, for higher-dimensional tangles and
higher-dimensional category theory.   Thus we begin with 
a brief sketch of both, concentrating on the situation at hand.  
We omit many details, which can be found in the references.  

\section{Preliminaries}

What are 2-tangles in 4 dimensions?  Roughly speaking, they are surfaces
in 4 dimensions going from one tangle in 3 dimensions to another tangle
in 3 dimensions.  Many ways of visualizing and representing them can be
found in the work of Carter, Rieger, and Saito \cite{CRS} and the
references therein.  To describe them more precisely, we equip the space
$[0,1]^4$ with standard coordinates $(x,y,z,t)$.  A 2-tangle $\alpha
\maps f \tto g$ in 4 dimensions is a compact 2-dimensional surface,
smoothly embedded in $[0,1]^4$ in such a way that its intersections with
the hyperplanes $\{t = 0\}$ and $\{t= 1\}$ are the tangles $f$ and $g$
in 3 dimensions.  The tangles $f,g$ may only touch the boundary of the
cube $[0,1]^3$ at the top and bottom, i.e., at the planes $\{z = 0\}$
and $\{z = 1\}$.  Near these planes, $f$ and $g$ must have a product
structure.   In other words, they must look like  straight vertical
lines.  Similarly, the surface $\alpha$ may only touch the  boundary of
$[0,1]^4$ at the hyperplanes $\{t = 0\}$, $\{t = 1\}$, $\{z = 0\}$, and
$\{z = 1\}$, and near these hyperplanes it must have a product
structure.  The surface $\alpha$ can have a boundary, but this boundary
must lie in the hyperplanes where either $t$ or $z$ is $0$ or $1$.  The
surface can also have right-angled corners, but only at those points
where both $t$ and $z$ are either $0$ or $1$.  Finally, we assume
without loss of generality that $\alpha$ is in `general position' in a
certain precise sense.  

Perhaps the most interesting examples of 2-tangles are compact surfaces
without boundary smoothly embedded in the interior of $[0,1]^4$.  
Topologists call these `knotted surfaces'.  The advantage of working
with more general 2-tangles is that it allows us to build complicated
knotted surfaces by gluing together simple 2-tangles.  We can glue
together 2-tangles in two basic ways.  First, given 2-tangles $\alpha
\maps f \tto g$ and $\beta \maps g \tto h$, we can form the 2-tangle
$\alpha\cdot  \beta \maps f \tto h$ by gluing together $\alpha$ and
$\beta$ along a hyperplane of constant $t$, namely, the hyperplane $t =
1$ for $\alpha$ and the hyperplane $t = 0$ for $\beta$.  Second, given
2-tangles $\alpha \maps f \tto g$ and $\beta \maps f' \tto g'$ such that
the composite tangles $ff'$ and $gg'$ are well-defined, we can form the
2-tangle $\alpha\circ\beta \maps ff' \tto gg'$ by gluing $\alpha$ and $\beta$
along a hyperplane of constant $z$.  

With a little work, we obtain an algebraic structure known as a
2-category, having objects, morphisms between objects, and 2-morphisms 
between morphisms.  We denote the 2-category of 2-tangles in 4 dimensions 
by $\T$.  The objects of $\T$ are certain equivalence classes of finite 
sets of points with distinct $y$ coordinates in the square $\{0 \le x,y \le
1\}$.  The morphisms in $\T$ are certain equivalence classes
of tangles, where it is crucial that tangles differing by a Reidemeister
move or by changing the relative heights of crossings, maxima, and
minima are not regarded as equivalent.  The 2-morphisms in $\T$ are
suitably defined ambient isotopy classes of 2-tangles.  

In any 2-category one can compose morphisms $f \maps A \to B$ and $g \maps B
\to C$ to obtain a morphism $fg \maps A \to C$.  (Note our ordering convention
here.)  In $\T$ this operation corresponds to composition of tangles.  
Also in any 2-category there are two ways to compose 
2-morphisms $\alpha$ and $\beta$, written $\alpha \cdot \beta$ and 
$\alpha\circ\beta$.  In $\T$ these operations work as described above.  
In the 2-categorical literature $\alpha\cdot\beta$ is usually 
called the `vertical composite' of $\alpha$ and $\beta$, while 
$\alpha\circ\beta$ is called the `horizontal composite'.  In what
follows, we write $1_f \circ \alpha$ simply as $f\alpha$, and
$\alpha \circ 1_f$ as $\alpha f$.  

Of course, there is a list of axioms to check \cite{KS} in order to 
show that $\T$ is a 2-category.  Actually, Kharlamov 
and Turaev \cite{KT} have already constructed a 2-category of
2-tangles in 4 dimensions.  For technical reasons ours is not 
quite the same, but we expect that it is `equivalent' in the sense 
of 2-category theory.

The reader may wonder why we did not consider a third basic way to
glue together 2-tangles: namely, along a hyperplane of constant $y$.
In fact, while the details are rather technical, this form of gluing
equips $\T$ with the the structure of a semistrict monoidal 2-category.
In a semistrict monoidal 2-category, one can tensor objects with
objects, morphisms or 2-morphisms, and there is an object $I$ serving
as the unit for the tensor product.  The reason for the term
`semistrict' is that certain equations which held strictly in a
monoidal category are weakened to 2-isomorphisms.  Most importantly,
in a monoidal category the equation $(A \tensor g)(f \tensor B') = (f
\tensor B)(A' \tensor g)$ holds for any morphisms $f \maps A \to A'$
and $g \maps B \to B'$, while in a semistrict monoidal 2-category this
is replaced by a specified 2-isomorphism, the `tensorator':
\[             {\bigotimes}_{f,g} \maps (A \tensor g)(f \tensor B') \tto
(f \tensor B)(A' \tensor g) . \]
Again there are various axioms that must hold.  These were first
explicitly listed by Kapranov and Voevodsky \cite{KV}, and later
expressed more tersely in the language of 2-category theory
\cite{BN,DS}.  While the details are rather lengthy, we can equip
$\T$ with unit object, tensor products, and tensorator, and check
that $\T$ becomes a semistrict monoidal 2-category.  (As noted
by Kharlamov and Turaev, Fischer's paper on 2-tangles \cite{Fischer} 
has serious flaws, such as not discussing the tensorator.)  

One can also consider gluing 2-tangles together 
along a hyperplane of constant $x$.  By the Eckmann-Hilton argument
\cite{BD,BN}, this form of gluing makes 
$\T$ into a semistrict braided monoidal 2-category.  This means, first
of all, that for any objects $A$ and $B$ there is a morphism $R_{A,B}
\maps A \tensor B \to B \tensor A$, called the `braiding'.  This
morphism must be a invertible up to a 2-isomorphism.  In addition, for
any morphisms $f \maps A \to A'$ and $g \maps B \to B'$, there are
braiding 2-isomorphisms
\[   R_{f,B} \maps (f \tensor B)R_{A',B} \tto R_{A,B}(B \tensor f) \]
and 
\[   R_{A,g} \maps (A \tensor g)R_{A,B'} \tto R_{A,B}(g \tensor A). \]
Now, in a strict braided monoidal category, for any objects $A,B$, and $C$
we have $(R_{A,B} \tensor C)(B \tensor R_{A,C}) = R_{A,B \tensor C}$ and
$(A \tensor R_{B,C})(R_{A,C} \tensor B) = R_{A \tensor B, C}$.  In fact,
these equations are crucial for proving the Yang-Baxter equation.  
In a semistrict braided monoidal 2-category these equations are weakened to 
specified 2-isomorphisms
\[  \tilde{R}_{(A|B,C)} \maps (R_{A,B} \tensor C)(B \tensor R_{A,C}) \tto 
R_{A,B \tensor C}  \]
and
\[    \tilde{R}_{(A,B|C)} \maps (A \tensor R_{B,C})(R_{A,C} \tensor B)
\tto R_{A \tensor B, C} .\]
There is also a list of axioms that must hold.  The first definition of
braided monoidal 2-category was given by Kapranov and Voevodsky 
\cite{KV}.  Later this definition was modified in various
ways by Baez and Neuchl \cite{BN}.   These modifications are necessary for 
the proper treatment
of 2-tangles, and especially for an unambiguous statement of the
Zamolodchikov tetrahedron equation, as had been noted by Breen
\cite{Breen}.  Subsequently Day and Street \cite{DS} re-expressed the
Baez-Neuchl definition in a more compact way, and Crans \cite{C} added
some axioms governing the braiding of the unit object.   In what
follows we use the definition given by Crans.   It turns
out that one can equip $\T$ with braiding morphisms and 2-isomorphisms
and check that it is a semistrict braided monoidal 2-category in this
sense.
 
There is a very special object $Z$ in $\T$, corresponding to a {\it
single point} in the square $\{0 \le x,y \le 1\}$.  Our result makes
precise the sense in which $\T$ is freely generated by this object.
For this we need to use the duality structure of $\T$.  The study of
duality in $n$-categories is only beginning, so before stating our
result we need to define braided monoidal 2-categories `with duals'.
Since we are working with unframed unoriented tangles, the object $Z$
has special properties: it is `self-dual' and `unframed'.  The concept
of a self-dual object is straightforward, but the concept of an
`unframed' object is rather subtle.  In the study of framed tangles, a
twist in the framing is often represented as a certain morphism $b_A
\maps A \to A$ known as the `balancing'.  In our situation, an
`unframed object' is not an object for which the balancing is the
identity, but one for which the balancing is {\it isomorphic} to the
identity via a 2-isomorphism that satisfies a highly nontrivial
equation of its own.  

Finally, the study of universal properties for $n$-categories is also
just beginning, so we must clarify what is meant by the `free' braided
monoidal 2-category with duals on one unframed self-dual object.  We
do so by means of a universal property.

In what follows, we take for granted our basic result that there
exists a semistrict braided monoidal 2-category $\T$ whose 2-morphisms
are ambient isotopy classes of smooth 2-tangles in 4 dimensions.  In
particular, letting $1_I$ denote the identity morphism of the unit
object of $\T$, the 2-morphisms $\alpha \maps 1_I \tto 1_I$ in $\T$
are precisely the ambient isotopy classes of compact surfaces without
boundary smoothly embedded in $[0,1]^4$, or equivalently, in $\R^4$.
In what follows we give an algebraic characterization of $\T$, and
thus of these `knotted surfaces'.

\section{Statement of Theorem}

In what follows, by monoidal and braided monoidal 2-categories we mean
`semistrict' ones as defined in reference \cite{BN}, but with the braided
monoidal 2-categories satisfying the extra axioms introduced
by Crans \cite{C}, which say that $R_{\cdot,\cdot}$, 
$\tilde R_{(\cdot|\cdot,\cdot)}$, and $\tilde R_{(\cdot,\cdot|\cdot)}$
are the identity whenever one of the arguments is the unit object $I$.

\begin{defn}\et A {\rm monoidal 2-category with duals}
is, to begin with, a monoidal 2-category equipped with the following
structures:
\begin{enumerate}
{\rm 
\item For every 2-morphism $\alpha \maps f \tto g$ there is a 
2-morphism $\alpha^*\maps g \tto f$ called the {\it dual} of $\alpha$.

\item For every morphism $f \from A\to B$ there is a morphism
$f^*\from B\to A$ called the {\it dual} of $f$, and 2-morphisms
$i_f\maps 1_A\tto f f^*$ and $e_f\maps f^* f\tto 1_B$,
called the {\it unit} and {\it counit} of $f$, respectively.

\item For any object $A$, there is a object $A^*$ called the {\it dual}
of $A$, morphisms $i_A\maps I\to A\tensor A^*$ and $e_A\maps
A^*\tensor A\to I$ called the {\it unit} and {\it counit} of $A$,
respectively, and a 2-morphism $T_A\maps (i_A\tensor
A)(A\tensor e_A)\tto 1_A$ called the {\it triangulator} of $A$.
}
\end{enumerate}

\noindent We say that a 2-morphism $\alpha$ is {\rm unitary} if it is
invertible and $\alpha^{-1} = \alpha^\ast$.
Given a 2-morphism $\alpha \maps f \tto g$, we define the {\rm 
adjoint} $\alpha^\dagger \maps g^* \tto f^*$ by 
\[  \alpha^\dagger= (g^*i_f) \cdot (g^*\alpha f^*) \cdot (e_g f^*) .\]

\noindent In addition, the structures above are also required to satisfy
the following conditions:
{\rm 
\begin{enumerate}

\item  $X^{**} = X$ for any object, morphism or 2-morphism $X$.

\item $1_X^* = 1_X$ for any object or morphism $X$.

\item  For all objects $A,B$, morphisms $f,g$, and 2-morphisms $\alpha,\beta$
for which both sides of the following equations are well-defined, we have
\[   (\alpha \cdot \beta)^* = \beta^*\cdot \alpha^*,\]
\[  (\alpha \hcomp \beta)^* = \alpha^*\hcomp \beta^*,\]
\[  (fg)^* = g^*f^*,\]
\[     (A\tensor \alpha)^* = A\tensor \alpha^*,  \qquad
(\alpha\tensor A)^* = \alpha^* \tensor A,  \]
\[    (A\tensor f)^* = A\tensor f^*, \qquad
 (f\tensor A)^* =  f^* \tensor A,  \]
and 
\[   (A\tensor B)^* = B^*\tensor A^* .\]

\item For all morphisms $f$ and $g$, the 2-morphism $\btensor_{f,g}$ is 
unitary. 

\item For any object or morphism $X$ we have
$i_{X^*} = e^*_X$ and $e_{X^*} = i^*_X$.  

\item For any object $A$, the 2-morphism $T_A$ is unitary.

\item If $I$ is the unit object, $T_I = 1_{1_I}$.

\item For any objects $A$ and $B$ we have 
\[   i_{A\tensor B} = i_A (A\tensor i_B \tensor A^*), \]
\[   e_{A\tensor B} = (B^*\tensor e_A \tensor B) e_B, \]
and
\[
T_{A\tensor B} = [(i_A\tensor A\tensor B) (A\tensor {\btensor}^{-1}_{i_B,e_A} 
\tensor B)(A \tensor B\tensor e_B)]
\cdot [(T_A \tensor B)\hcomp (A\tensor T_B)] .  \]

\item For any object $A$ and morphism $f$ we have
\[      i_{A \tensor f} = A \tensor i_f , \qquad
     i_{f \tensor A} = i_f \tensor A ,\]
\[      e_{A \tensor f} = A \tensor e_f , \qquad
     e_{f \tensor A} = e_f \tensor A .\] 
  
\item For any morphisms $f$ and $g$, $i_{fg} = i_f\cdot (f i_g
f^*)$ and $e_{fg} = (g^* e_f g) \cdot e_g$.

\item For any morphism $f$, $i_f f\cdot fe_f = 1_f$ and $f^* i_f\cdot 
e_f f^* = 1_{f^*}$.

\item For any 2-morphism $\alpha$, $\alpha^{\dagger *} = 
\alpha^{* \dagger}$.

\item For any object $A$ we have
\[ [i_A(A\tensor T^{\dagger }_{A^*})] \cdot
[{\btensor}_{i_A,i_A}^{-1}(A\tensor e_A\tensor A^*)]\cdot 
[i_A(T_A \tensor A^*)] = 1_{i_A}.\]
\end{enumerate}
}
\end{defn}

The final equation has the following geometrical interpretation in
terms of movie moves (see \cite{CRS}):

\vbox{
\medskip
\centerline{\epsfysize=2in\epsfbox{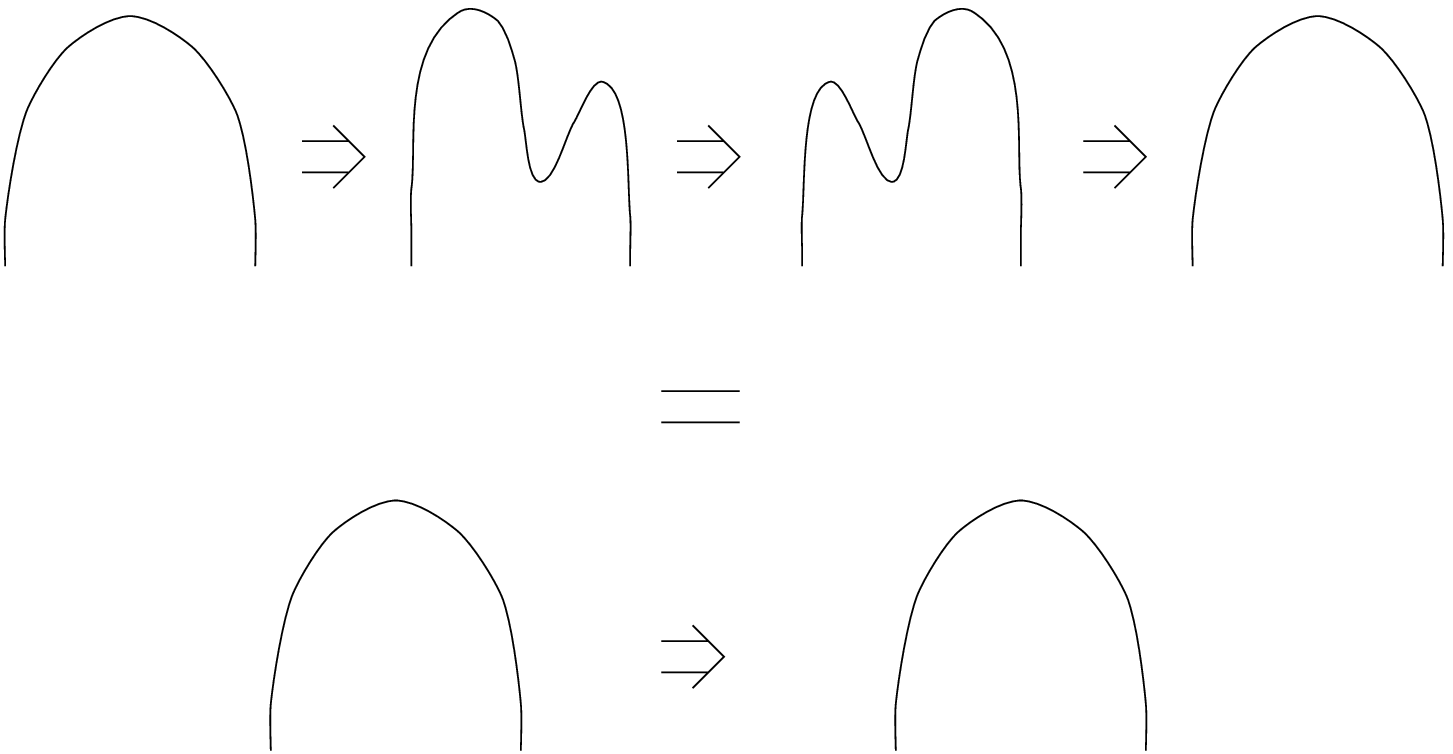}}
\bigskip
\centerline{1.  Equation satisfied by the triangulator}
\medskip
}

\noindent Note that the equations in clause 5 of the definition above
allow us to express the counit in terms of the unit (or vice versa)
using duality.  This allows us to avoid mentioning the counit at certain
points in some definitions below.

\begin{defn} \label{bm2cat.w.duals}
\et A {\rm braided monoidal 2-category with duals} is a
monoidal 2-category with duals that is also a braided monoidal 2-category
for which the braiding is unitary in the sense that: 
\begin{enumerate} 
\item For any objects $A,B$, the 2-morphisms
$i_{R_{A,B}}$ and $e_{R_{A,B}}$ are unitary.  
\item For any object $A$ and morphism $f$, 
the 2-morphisms $R_{A,f}$ and $R_{f,A}$ are unitary.
\item For any objects $A,B,C$, the 2-morphisms
$\tilde{R}_{(A,B|C)}$ and
$\tilde{R}_{(A|B,C)}$ are
unitary.
\end{enumerate}
\end{defn}

In general the above definition does not deal adequately with the
subtle issue of framings (or algebraically speaking, balancings).  We
can use this definition here because we are mainly interested in $\T$,
which is generated by an `unframed self-dual' object in the sense
defined below.  Geometrically, this object is simply a point embedded
in the unit square.

In previous work \cite{B}, it was shown how the balancing 
arises naturally in any braided monoidal category with duals.
The same idea applies to braided monoidal 2-categories with duals.
Explicitly, for any object $A$ in a braided monoidal 2-category with
duals, the balancing $b_A \maps A \to A$ is given by:
\[  b_A = (e_A^* \tensor A)(A^* \tensor R_{A,A})(e_A \tensor A) .\]
For an `unframed' object $A$, the 1st Reidemeister move corresponds
to a 2-morphism $V_A \maps b_A \tto 1_A$.  
However, the connection to the movie moves of Carter, 
Rieger and Saito becomes a bit clearer if we work not with the balancing but
with the closely related morphism $i_{A^*}R_{A^*,A} \maps 1 \to A \tensor A^*$.
The 1st Reidemeister move then corresponds to a 2-isomorphism 
\[           W_A \maps i_A \tto  i_{A^*}R_{A^*,A} \]
which we call the `writhing', with the following geometrical interpretation:

\vbox{
\medskip
\centerline{\epsfysize=0.75in\epsfbox{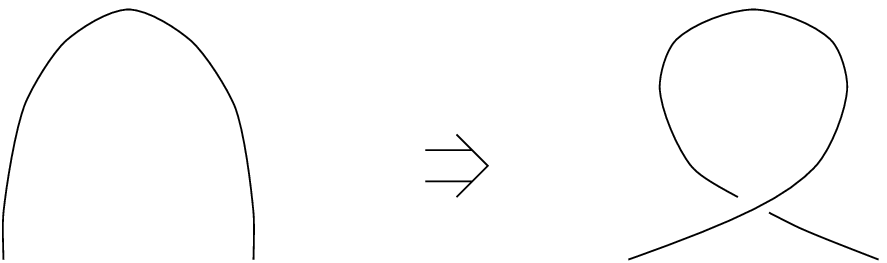}}
\bigskip
\centerline{2.  The writhing}
\medskip
}

\noindent One can construct a 2-isomorphism $V_A \maps b_A \tto 1_A$ given a
2-isomorphism $W_A \maps i_A \tto i_{A^*}R_{A^*,A}$, and conversely.
In what follows we only study the writhing for a self-dual object $A$.  

\begin{defn} \et A {\rm self-dual object} 
in a braided monoidal 2-category with duals is an object $A$ with
$A^* = A$. \end{defn}

\begin{defn}\et A self-dual object $A$
in a braided monoidal 2-category with duals is {\rm unframed} 
if $\tilde R_{(A|A,A)} = 1$, $\tilde R_{(A,A|A)} = 1$, and it
is equipped with a unitary 2-morphism
\[  W_A \maps i_A \tto i_{A} R_{A,A} \]
called the {\rm writhing}, satisfying the equation:
\[
T_A^\dagger \cdot 
((A\tensor W_A)(e_{A}\tensor A)) \cdot 
((A\tensor (i_{A}R_{A,A}))R^\dagger_{A,i_A})\cdot \]
\[ 
((A\tensor (i_{A}R_{A,A}))\tilde{R}^\dagger_{(A|A,A)}
(A\tensor  e_{A})) \cdot 
((A\tensor (i_{A}i^*_{R_{A,A}}))(R^*_{A,A}\tensor A)(A\tensor e_{A}))
\]
\[  =T_A^{-1} \cdot
((i_A\tensor A)(A\tensor (i_{R_{A,A}} e_A)))\cdot ((i_A\tensor A)
(A\tensor (R_{A,A} W_A^\dagger)))\cdot \]
\[({R^\dagger}^{-1}_{e_{A},A}(A\tensor (R_{A,A}e_{A})))\cdot 
((A\tensor i_A) \tilde{R}^\dagger_{(A,A|A)}
(A\tensor (R_{A,A}e_{A}))) \cdot 
\]
\[((A\tensor i_A)(R^*_{A,A}\tensor A)(A\tensor (e_{R_{A,A}}e_{A}))) \]
\end{defn}

The rather terrifying equation above has the following geometrical
interpretation in terms of `movie moves': 

\vbox{
\medskip
\centerline{\epsfysize=2in\epsfbox{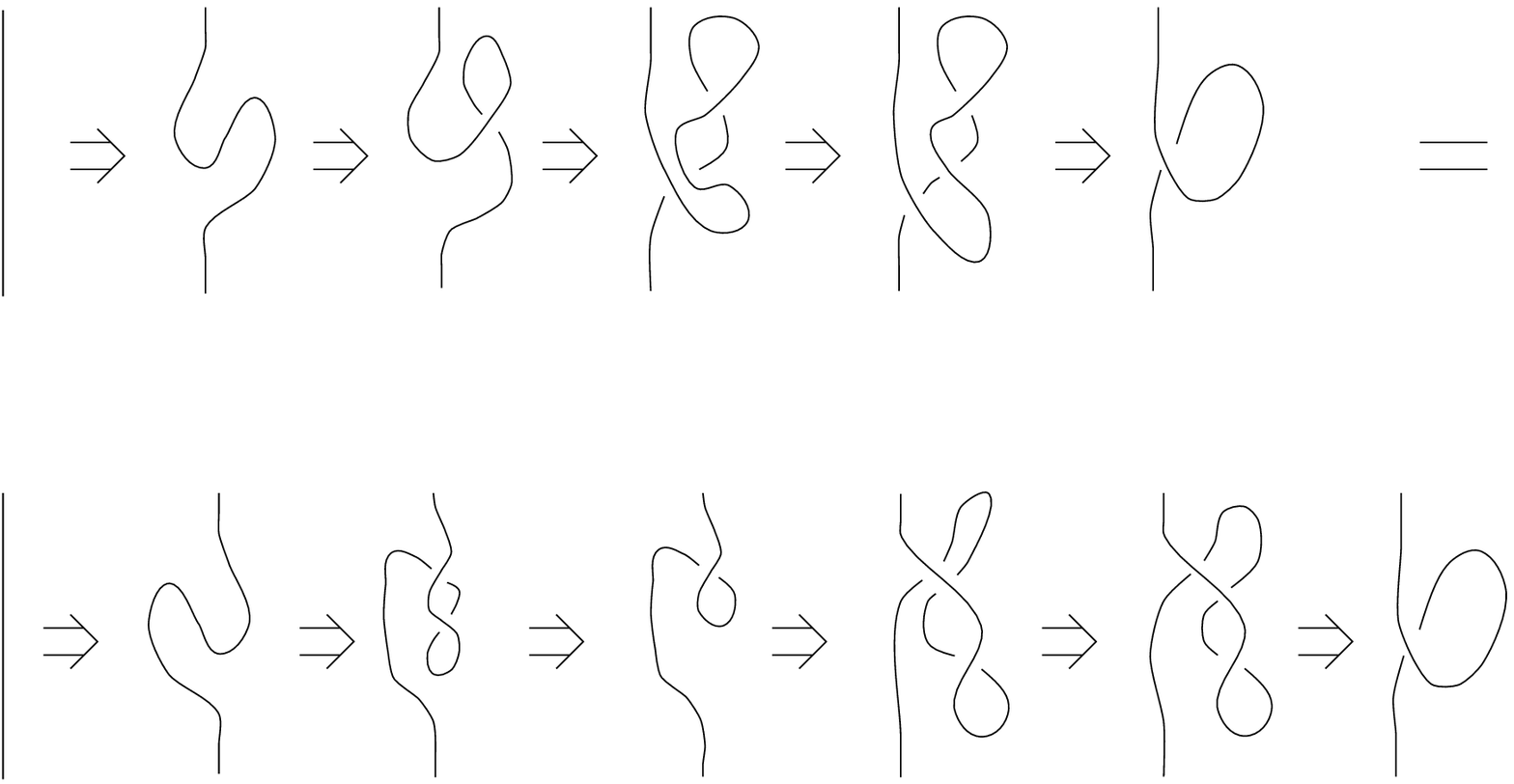}}
\bigskip
\centerline{3.  Equation satisfied by the writhing}
\medskip
}

\noindent 
The conditions $\tilde R_{(A|A,A)} = 1$ and
$\tilde R_{(A|A,A)} = 1$ have nothing to do with framing per se;
they actually amount to a kind of `strictification' of the braiding
as far as the object $A$ is concerned.  We include them in
the definition of `unframed object' merely to simplify the exposition
in what follows, and will remove them in our forthcoming more
detailed treatment.

\begin{defn}\et We say a braided monoidal 2-category is {\rm generated}
by an unframed self-dual object $Z$ if:

\begin{enumerate}
{\rm 
\item Every object is a tensor product of copies of $Z$.
\item Every morphism can be obtained by composition from:
\begin{alphalist}
\item $1_Z$, 
\item $i_Z$, 
\item $R_{Z,Z}$, 
\item tensor products of arbitrary objects with the above morphisms,
\item duals of the above morphisms.
\end{alphalist}
\item Every 2-morphism can be obtained by horizontal and vertical composition
from: 
\begin{alphalist}
\item 2-morphisms $1_f$ for arbitrary morphisms $f$,
\item 2-morphisms $\bigotimes_{f,g}$ for arbitrary morphisms $f$ and $g$, 
\item 2-morphisms $R_{Z,f}$ and $R_{f,Z}$ for arbitrary morphisms $f$, 
\item 2-morphisms $i_f$ for arbitrary morphisms $f$, 
\item $T_Z$, 
\item $W_Z$,
\item tensor products of arbitrary objects with the above 2-morphisms, and
\item duals of the above 2-morphisms.
\end{alphalist}
}
\end{enumerate}
\end{defn}

\begin{defn}\et For monoidal 2-categories $\C$ and $\D$, a {\rm strict 
monoidal 2-functor} $F\maps \C
\to \D$ is a 2-functor such that $F(1) = 1$,
$F(A\tensor X) = F(A)\tensor F(X)$,
$F(X\tensor A) = F(X)\tensor F(A)$ and $F(\btensor_{f,g}) =
\btensor_{F(f),F(g)}$, for any object $A$, object, morphism or 2-morphism
$X$, and morphisms $f$ and $g$.  
\end{defn}

\begin{defn} \et If $\C,\D$ are braided monoidal 2-categories with duals
and $\C$ is generated by the unframed self-dual object $Z$, we say a monoidal
2-functor $F \maps \C \to \D$ mapping $Z$ to an unframed self-dual object
in $\D$ {\rm preserves braiding and duals strictly on 
the generator} if: 
{\rm 
\begin{alphalist}
\item  $F(X^*) = F(X)^*$
for every object, morphism, or 2-morphism
$X$, 
\item $F(i_f) = i_{F(f)}$ for every morphism $f$,
\item $F(i_Z) = i_{F(Z)}$, 
\item $F(T_Z) = T_{F(Z)}$,
\item $F(W_Z) = W_{F(Z)}$
\item $F(R_{Z,Z}) = R_{F(Z),F(Z)}$, 
\item $F(R_{Z,f}) = R_{F(Z),F(f)}$ and $F(R_{f,Z}) = R_{F(f), F(Z)}$ for 
$f$ equal to $R_{Z,Z}, R^*_{Z,Z}, i_Z,$ and $e_Z$.  
\end{alphalist}
}
\end{defn} 

\begin{thm}\et
There is a braided monoidal 2-category with duals $\T$ for which there
is an explicit one-to-correpondence between 2-morphisms of $\T$ and
smooth ambient isotopy classes of unframed unoriented smooth 2-tangles
in 4 dimensions.  ${\cal T}$ is generated by an
unframed self-dual object $Z$, and
for any braided monoidal 2-category with duals $\C$ and unframed
self-dual object $A \in \C$ there is a unique strict monoidal 2-functor 
$F\maps \T \to \C$ with $F(Z) = A$ 
that preserves braiding and duals strictly on the generator.
\end{thm}

Thus we say that $\T$ is the {\it free
braided monoidal 2-category with duals on one unframed self-dual object}.

\section{Conclusions}

Our result implies that we obtain an invariant of unframed unoriented
2-tangles in 4 dimensions from any braided monoidal 2-category with
duals containing a self-dual unframed object.  Similarly, we expect to
be able to prove that the 2-category of framed oriented 2-tangles is
the free braided monoidal 2-category with duals on one object.  This would
enable us to obtain invariants of framed oriented 2-tangles from
arbitrary braided monoidal 2-categories with duals.  However, given the rather
lengthy definitions above, it is natural to wonder whether we can
actually find invariants this way in practice: are there any interesting
{\it examples} of braided monoidal 2-categories with duals?

We believe there are many examples and that the problem is mainly a
matter of developing the machinery to get our hands on them.  First,
there is plenty of evidence \cite{BD,BD2} suggesting that we can
obtain braided monoidal 2-categories from the homotopy 2-types of
double loop spaces.  Second, Neuchl and the first author have shown
how to obtain braided monoidal 2-categories from monoidal 2-categories
by a `quantum double' construction \cite{BN}.  It seems plausible that
applying this construction to a monoidal 2-category with duals will
give a braided monoidal 2-category with duals.  This reduces the
question to obtaining monoidal 2-categories with duals.  A good
example of one of these should be the monoidal 2-category of unitary 
representations of a 2-groupoid, just as the monoidal category of unitary 
representations of a groupoid is an example of a monoidal category 
with duals \cite{B}.  Third, Crane and Frenkel have sketched a way to 
construct Hopf categories from Kashiwara and Lusztig's canonical bases 
for quantum groups \cite{CF}.  There is reason to hope that the representation
2-categories of these Hopf categories are monoidal 2-categories with
duals.  Fourth, just as one can construct braided monoidal categories
from solutions of the Yang-Baxter equation, one can construct braided
monoidal 2-categories from solutions of the Zamolodchikov tetrahedron
equations \cite{KV}.  Many such solutions are known \cite{CS}, so one
may hope that some give braided monoidal 2-categories with duals.
Finally, one expects `braided monoidal 3-Hilbert spaces' to be
interesting examples of braided monoidal 2-categories with duals
\cite{B}.

\subsection*{Acknowledgements}

Many of the ideas underlying this work were developed with James Dolan.
We would also like to thank Scott Carter and Masahico Saito for a
great deal of useful correspondence.

\end{document}